# QUANTILE REGRESSION IN PARTIALLY LINEAR VARYING COEFFICIENT MODELS

By Huixia Judy Wang[1], Zhongyi Zhu[2] and Jianhui Zhou

*North Carolina State University, Fudan University and University of Virginia*

Semiparametric models are often considered for analyzing longitudinal data for a good balance between flexibility and parsimony. In this paper, we study a class of marginal partially linear quantile models with possibly varying coefficients. The functional coefficients are estimated by basis function approximations. The estimation procedure is easy to implement, and it requires no specification of the error distributions. The asymptotic properties of the proposed estimators are established for the varying coefficients as well as for the constant coefficients. We develop rank score tests for hypotheses on the coefficients, including the hypotheses on the constancy of a subset of the varying coefficients. Hypothesis testing of this type is theoretically challenging, as the dimensions of the parameter spaces under both the null and the alternative hypotheses are growing with the sample size. We assess the finite sample performance of the proposed method by Monte Carlo simulation studies, and demonstrate its value by the analysis of an AIDS data set, where the modeling of quantiles provides more comprehensive information than the usual least squares approach.

**1. Introduction.** Various nonparametric models have been developed for longitudinal data analysis. One popular nonparametric specification is the varying coefficient model, where the coefficients are smooth nonparametric functions of some factors such as measurement time. Varying coefficient models were proposed by Hastie and Tibshirani [8], and later extended to longitudinal studies by Chiang, Rice and Wu [3], Huang, Wu and Zhou

Received October 2008; revised December 2008.
[1]Supported by NSF Grant DMS-07-06963.
[2]Supported by National Natural Science Foundation of China Grant 10671038 and Shanghai Leading Academic Discipline Project B210.

[13], and Qu and Li [21], among others. Though flexible, the varying coefficient models may overfit data when some covariate effects are indeed time-invariant. This motivates the partially linear varying coefficient model (PLVC)

$$y = x'\alpha(t) + z'\beta + e, \tag{1.1}$$

where $\alpha(t)$ comprises $p$ unknown smooth functions, $\beta$ is a $q$-dimensional parameter vector, and $e$ is the random error.

We consider PLVC models for longitudinal data due to their good balance between flexibility and parsimony. The PLVC models have been studied by Ahmad, Leelahanon and Li [1] and Fan and Huang [4] for cross-sectional data, and by Sun and Wu [23] and Fan, Huang and Li [5] for longitudinal data. The current literature is mainly confined to estimating the conditional mean function of $y$. The focus of this paper is to estimate and conduct inference on the conditional quantile curves without any specification of the error distribution or intra-subject dependence structure.

Suppose that we have $n$ subjects, and the $i$th subject has $m_i \geq 1$ repeated measurements over time. At a given quantile level $\tau \in (0,1)$, we assume the following PLVC quantile regression model:

$$y_{ij} = x'_{ij}\alpha(t_{ij}, \tau) + z'_{ij}\beta(\tau) + e_{ij}(\tau), \qquad i = 1, \ldots, n, j = 1, \ldots, m_i, \tag{1.2}$$

where $y_{ij}$ denotes the $j$th outcome of the $i$th subject, $\alpha(t, \tau) = (\alpha_1(t, \tau), \ldots, \alpha_p(t, \tau))'$ are unknown smooth functions of $t$ for all $t \in R$, $x_{ij} = (x_{ij,1}, \ldots, x_{ij,p})' \in \mathbb{R}^p$ and $z_{ij} = (z_{ij,1}, \ldots, z_{ij,q})' \in \mathbb{R}^q$ are the design vectors for the time-varying coefficients and constant coefficients, respectively, and $e_{ij}(\tau)$ is the random error whose $\tau$th quantile conditional on $(x, z, t)$ equals zero. We assume that the observations, and therefore $e_{ij}$, are dependent within the same subjects, but independent across subjects. The form of the error distribution and the intra-subject dependence structure are left unspecified.

The quantile PLVC model is a valuable alternative to the conditional mean models for analyzing longitudinal data. First, fitting data at a set of quantiles provides a more comprehensive description of the response distribution than does the mean. In many applications, the functional impacts of the covariates on the response may vary at different percentiles of the distribution. For instance, the analysis of an AIDS data set in Section 6 shows that the effect of the initial measurement on CD4 percentage is time-decaying for severely ill people (in the lower tail of the CD4 distribution), whereas it tends to be constant in the upper quantiles. Such important features would be overlooked by the regression approaches that focus only on the mean or the median. Second, the modeling of different conditional quantiles can be used to construct prediction intervals; see, for instance, [2] and [19]. Third, when the center of the conditional distribution of $y$ is of interest, the median regression, a special case of quantile regression, provides more robust estimators than the mean regression. In addition, quantile regression does not assume any parametric form on the error distribution, and thus is able to accommodate nonnormal errors, as often seen in longitudinal studies.

Even though linear quantile regression has been well developed, theory and methodology for nonparametric and semiparametric quantile models are lagging. Hendricks and Koenker [11], Yu and Jones [26], among others, studied nonparametric quantile regression for independent observations. Recently, several authors studied quantile regression for varying coefficient models. Cai and Xu [2] considered local polynomial estimators for time series data. Honda [12] and Kim [15] studied varying coefficient models for independent data using local polynomials and splines, respectively.

The present paper is the first to develop theory and methodology for analyzing longitudinal data in the quantile PLVC models. Koenker [17] and Lipsitz et al. [18] discussed quantile regression for longitudinal data in linear models. Wei and He [25] studied semiparametric quantile autoregression models for longitudinal data, where the errors $e_{ij}$ are treated as independent in the hypothesis testing.

In this paper, we propose to estimate the quantile smooth coefficients using basis function approximations. Under some regularity conditions, we obtain the optimal convergence rate of the functional coefficient curves $\hat{\alpha}(t, \tau)$, and establish the asymptotic normality of $\hat{\beta}(\tau)$. Even though a Wald-type test can be constructed using the asymptotic distribution of $\hat{\beta}(\tau)$, its finite sample performance is very sensitive to the estimation of the error density function evaluated at the quantile of interest. To avoid this problem, we develop a quantile rank score test for inference on $\beta(\tau)$, and show that it is superior to the Wald-type test in finite samples. Focusing on a varying coefficient model, Kim [15] implemented a similar inference procedure for testing the hypothesis that all of the varying coefficients are constant. However, a question of more practical use is to test whether a particular covariate or a subset of covariates have time-varying effects. Extensions of the rank score test to this type of hypothesis testing problems are theoretically challenging, because the dimensions of the parameter spaces under both the null and the alternative hypotheses are increasing with the sample size. We provide the asymptotic results under this practically relevant setting, making it possible to test the constancy of the coefficients one at a time, a necessary step in selecting varying-coefficient models of different complexities.

In Section 2, we present the proposed estimation procedure and the large sample properties of the resulting estimators. We develop inferential procedures for testing $\beta(\tau)$ in Section 3, and for testing the constancy of a subset of $\alpha_l(t, \tau)$'s in Section 4. We assess the finite sample performance of the proposed procedures with simulation studies in Section 5. The merit of the proposed methods is illustrated by analyzing a CD4 depletion data in Section 6. Section 7 concludes the paper with some discussion. Proofs are deferred to the Appendix.

## 2. The proposed method.

2.1. *Estimation.* For the ease of presentation, we will omit $\tau$ from $\alpha_l(t,\tau)$, $\beta(\tau)$ and $e_{ij}(\tau)$ in model (1.2) wherever clear from the context, but we should bear in mind that those quantities are $\tau$-specific. Without loss of generality, we assume that $t_{ij} \in [0,1]$ for all $i$ and $j$ throughout.

Let $\pi(t) = (B_1(t), \ldots, B_{k_n+\hbar+1}(t))'$ be a set of B-spline basis functions of order $\hbar + 1$ with $k_n$ quasi-uniform internal knots. We approximate each $\alpha_l(t)$ by a linear combination of normalized B-spline basis functions $\alpha_l(t) \approx \sum_{s=1}^{k_n+\hbar+1} B_s(t)\theta_{l,s} = \pi(t)'\theta_l$, where $\theta_l = (\theta_{l,1}, \ldots, \theta_{l,k_n+\hbar+1})'$ is the spline coefficient vector. For details on the construction of B-spline basis functions, the readers are referred to Schumaker [22]. With the B-spline basis, model (1.2) can be approximated by

$$y_{ij} \approx \sum_{l=1}^{p} \sum_{s=1}^{k_n+\hbar+1} x_{ij,l} B_s(t_{ij}) \theta_{l,s} + \sum_{s=1}^{q} z_{ij,s} \beta_s + e_{ij} = \Pi'_{ij}\Theta + z'_{ij}\beta + e_{ij},$$

where $\Pi_{ij} = (x_{ij,1}\pi'_{ij}, \ldots, x_{ij,p}\pi'_{ij})' \in \mathbb{R}^{p_{k_n}}$, $\Theta = (\theta_{l,s}) \in \mathbb{R}^{p_{k_n}}$ with $p_{k_n} = p(k_n + \hbar + 1)$, and $\pi_{ij} = \pi(t_{ij})$. The quantile coefficient estimates $\hat{\Theta}$ and $\hat{\beta}$ can be obtained by minimizing

$$(2.1) \qquad \sum_{i=1}^{n} \sum_{j=1}^{m_i} \rho_\tau(y_{ij} - \Pi'_{ij}\Theta - z'_{ij}\beta),$$

where $\rho_\tau(u) = u\{\tau - I(u < 0)\}$ is the quantile loss function. This estimation method is a one-step procedure applicable to cases where the measurements are either regularly or irregularly observed.

In practice, we choose lower order splines, such as $\hbar = 2$ or $\hbar = 3$ corresponding to quadratic or cubic splines. For demonstration, we assume that the number of internal knots, $k_n$, is the same for each varying coefficient, even though it is allowed to vary in real applications; see Section 5.1 for a model selection criterion for determining $k_n$.

2.2. *Large sample properties.* Before presenting the main asymptotic results, we first introduce two definitions.

DEFINITION 1. Define $\mathcal{H}_r$ as the collection of all functions on $[0,1]$ whose $m$th order derivative satisfies the Hölder condition of order $\nu$ with $r \equiv m+\nu$. That is, for any $h \in \mathcal{H}_r$, there exists a constant $c \in (0, \infty)$ such that for each $h \in \mathcal{H}_r$, $|h^{(m)}(s) - h^{(m)}(t)| \leq c|s-t|^\nu$, for any $0 \leq s, t \leq 1$.

DEFINITION 2. The function $g(x,t)$ is said to belong to the varying coefficient class of functions $\mathcal{Y}$ if (i) $g(x,t) = x'h(t) \equiv \sum_{l=1}^{p} x_l h_l(t)$; (ii) $\sum_{l=1}^{p} E\{x_l \times h_l(t)\}^2 < \infty$, where $x_l$ and $h_l(t) \in \mathcal{H}_r$ are the $l$th coordinates of $x$ and $h(t)$, respectively, $l = 1, \ldots, p$.

For convenience, we define $Z_i = (z_{i1}, \ldots, z_{im_i})'$ as the $m_i \times q$ design matrix on the $i$th subject, and $Z = (Z_1', \ldots, Z_n')'$. Similarly, denote $X_i = (x_{i1}, \ldots, x_{im_i})'$, $T_i = (t_{i1}, \ldots, t_{im_i})' \in \mathbb{R}^{m_i}$, $X = (X_1', \ldots, X_n')'$ and $T = (T_1', \ldots, T_n')'$. Let $Z_{i,l}$ be the $l$th column of $Z_i$, $F_{ij}$ be the cumulative distribution function (CDF) and $f_{ij}$ be the density function of $e_{ij}$ conditional on $(x_{ij}, z_{ij}, t_{ij})$. We denote $B_i = \mathrm{diag}(f_{i1}(0), \ldots, f_{im_i}(0))$ and $B = \mathrm{diag}(B_1, \ldots, B_n)$.

To obtain the asymptotic distribution of $\hat{\beta}$, we first need to adjust for the dependence of $Z$ and $(X,T)$, which is a common complication in semiparametric models. Similar to [1], we denote $\phi(x_{ij}, t_{ij}) = \sum_{l=1}^{p} x_{ij,l} h_l(t_{ij}) \in \mathcal{Y}$ and $\phi(X_i, T_i) = (\phi(x_{i1}, t_{i1}), \ldots, \phi(x_{im_i}, t_{im_i}))'$. Let

$$(2.2) \quad \phi_l^*(\cdot, \cdot) = \arg\inf_{\phi \in \mathcal{Y}} \sum_{i=1}^{n} E[\{Z_{i,l} - \phi(X_i, T_i)\}' B_i \{Z_{i,l} - \phi(X_i, T_i)\}]$$

and $m_l(x,t) = E(Z_{i,l} | X = x, T = t)$, $l = 1, \ldots, q$. Note that

$$\sum_{i=1}^{n} E[\{Z_{i,l} - \phi(X_i, T_i)\}' B_i \{Z_{i,l} - \phi(X_i, T_i)\}]$$

$$= \sum_{i=1}^{n} E[\{Z_{i,l} - m_l(X_i, T_i)\}' B_i \{Z_{i,l} - m_l(X_i, T_i)\}]$$

$$+ \sum_{i=1}^{n} E[\{m_l(X_i, T_i) - \phi(X_i, T_i)\}' B_i \{m_l(X_i, T_i) - \phi(X_i, T_i)\}].$$

Therefore, $\phi_l^*(X_i, T_i)$ are the projections of $m_l(X_i, T_i)$ onto the varying coefficient functional space $\mathcal{Y}$ (under the $L_2$-norm). In other words, $\phi_l^*(X_i, T_i)$ is an element that belongs to $\mathcal{Y}$ and it is the closest function to $m_l(X_i, T_i)$ among all the functions in $\mathcal{Y}$, for any $l = 1, \ldots, q$. In (2.2), we consider the weighted projection to account for the heteroscedasticity through $B_i$.

We define

$$K_n = \sum_i \{Z_i - \phi^*(X_i, T_i)\}' B_i \{Z_i - \phi^*(X_i, T_i)\},$$

$$\Lambda_n = \sum_i \{Z_i - \phi^*(X_i, T_i)\}' A_i(\Delta) \{Z_i - \phi^*(X_i, T_i)\},$$

$$A_i(\Delta) = \mathrm{Cov}(\psi_\tau(e_{i1}), \ldots, \psi_\tau(e_{im_i})),$$

where $\phi^*(X_i, T_i) = (\phi_1^*(X_i, T_i), \ldots, \phi_q^*(X_i, T_i))$ and $\psi_\tau(u) = \tau - I(u < 0)$. The $A_i(\Delta)$ is an $m_i \times m_i$ symmetric matrix with the $(j_1, j_2)$-element $\tau - \tau^2$ if $j_1 = j_2$, and $\delta_{ij_1j_2} - \tau^2$ if $j_1 \neq j_2$, where $\delta_{ij_1j_2} = P(e_{ij_1}(\tau) < 0, e_{ij_2}(\tau) < 0)$ measures the tail dependence of a pair of residuals from the same subject, and $\Delta$ is the collection of all $\delta_{ij_1j_2}$'s. The following assumptions are needed to obtain the asymptotic properties of $\hat{\alpha}(t)$ and $\hat{\beta}$.

A1 For some $r \geq 1$, $\alpha_l(t) \in \mathcal{H}_r, l = 1, \ldots, p$.
A2 The conditional distribution of $T$ given $X = x$ has a bounded density of $f_{T|X}$ satisfying $0 < c_1 \leq f_{T|X}(t|x) \leq c_2 < \infty$, uniformly in $x$ and $t$ for some positive constants $c_1$ and $c_2$.
A3 The numbers of measurements $m_i$ are uniformly bounded for all $i = 1, \ldots, n$.
A4 Uniformly over $i$ and $j$, $f_{ij}(\cdot)$ is bounded from infinity, and it is bounded away from zero and has a bounded first derivative in the neighborhood of zero.
A5 For all $i$ and $j$, the random design vectors $x_{ij}, z_{ij}$ are bounded in probability.
A6 Let $N = \sum_{i=1}^n m_i$ denote the total number of observations. The eigenvalues of $N^{-1}K_n$ and $N^{-1}\Lambda_n$ are bounded away from infinity and zero for sufficiently large $n$.

THEOREM 1. *Under assumptions A1–A6, if $r \geq 1$ and $k_n \approx n^{1/(2r+1)}$, then*

$$(2.3) \quad \frac{1}{N} \sum_{i=1}^n \sum_{j=1}^{m_i} \{\hat{\alpha}_l(t_{ij}) - \alpha_l(t_{ij})\}^2 = O_p(n^{-2r/(2r+1)}), \quad l = 1, \ldots, p,$$

*and*

$$(2.4) \quad \Lambda_n^{-1/2} K_n(\hat{\beta} - \beta_0) \to N(0, I_q).$$

*Throughout, we use $a_n \approx b_n$ to mean that $a_n$ and $b_n$ have the same order as $n \to \infty$.*

Assumption A1 is required to achieve the optimal convergence rate of $\hat{\alpha}_l(t)$. Assumption A2 ensures that $k_n \Pi' B \Pi$ is positive definite for sufficiently large $n$, which is needed for proving the consistency of the estimators. Assumptions A3 and A4 are standard assumptions used in longitudinal studies and quantile regression, respectively. The boundness assumption in A5 is made for convenience. It suffices to assume that $x_{ij}$ and $z_{ij}$ have bounded fourth moments, but this will complicate the technical proof. Assumption A6 is used to represent the asymptotic covariance matrix of $\hat{\beta}$ and to obtain the optimal convergence rate of the estimators of the parametric and nonparametric parts.

**3. Inference on nonvarying coefficients.** In this section, we propose two large sample inference procedures for testing the nonvarying coefficients $\beta$, including a Wald-type test and a rank-score-based test.

3.1. *Wald-type test.* Based on the asymptotic normality (2.4), a Wald-type test can be constructed for inference on $\beta$ through direct estimation of the covariance matrix, which involves the unknown error density $f_{ij}(0)$. Koenker [16] discussed several ways to estimate $f_{ij}(0)$ and provided recommendations for selecting the smoothing parameters. In this paper, we adopt the idea of Hendricks and Koenker [11] and estimate $f_{ij}(0)$ by the difference quotient,

$$
\begin{aligned}
(3.1) \quad \hat{f}_{ij}(0) = 2\epsilon_n [x'_{ij}\{\hat{\alpha}(t_{ij}, \tau + \epsilon_n) - \hat{\alpha}(t_{ij}, \tau - \epsilon_n)\} \\
+ z'_{ij}\{\hat{\beta}(\tau + \epsilon_n) - \hat{\beta}(\tau - \epsilon_n)\}]^{-1},
\end{aligned}
$$

where $\epsilon_n$ is a bandwidth parameter tending to zero as $n \to \infty$. Throughout our numerical studies, we choose $\epsilon_n = 1.57 n^{-1/3}(1.5\phi^2\{\Phi^{-1}(\tau)\}/[2\{\Phi^{-1}(\tau)\}^2 + 1])^{2/3}$ following Hall and Sheather [7], where $\Phi(\cdot)$ and $\phi(\cdot)$ are the CDF and density function of the standard normal distribution.

Denote

$$
\begin{aligned}
&\hat{B}_i = \text{diag}\{\hat{f}_{i1}(0), \ldots, \hat{f}_{im_i}(0)\}, \qquad \hat{B} = \text{diag}(\hat{B}_1, \ldots, \hat{B}_n), \\
(3.2) \quad &\Pi = (\Pi_{11}, \ldots, \Pi_{1m_1}, \ldots, \Pi_{nm_n})', \\
&\hat{P} = \Pi(\Pi'\hat{B}\Pi)^{-1}\Pi'\hat{B}, \qquad \hat{Z}^* = (I - \hat{P})Z,
\end{aligned}
$$

and let $\hat{Z}_i^*$ be the rows of $\hat{Z}^*$ corresponding to the $i$th subject.

THEOREM 2. *Let $\hat{e}_{ij} = y_{ij} - x'_{ij}\hat{\alpha}(t_{ij}) - z'_{ij}\hat{\beta}$, $\hat{e}_i = (\hat{e}_{i1}, \ldots, \hat{e}_{im_i})'$,*

$$
(3.3) \quad \hat{\Lambda}_n = \sum_{i=1}^n \hat{Z}_i^{*\prime} \psi(\hat{e}_i)\psi'(\hat{e}_i)\hat{Z}_i^*, \qquad \hat{K}_n = \sum_{i=1}^n \hat{Z}_i^{*\prime}\hat{B}_i\hat{Z}_i^*.
$$

*If $\epsilon_n \to 0$, $\liminf_{n \to \infty} n\epsilon_n^2 > 0$ and the assumptions of Theorem 1 hold, then $N^{-1}(\hat{K}_n - K_n) = o_p(1)$ and $N^{-1}(\hat{\Lambda}_n - \Lambda_n) = o_p(1)$.*

3.2. *Rank score test.* Theorem 2 shows that the covariance matrix of $\hat{\beta}$ can be estimated consistently. However, the finite sample performance of the Wald-type test is sensitive to the choice of $\epsilon_n$. As an alternative, we consider a quantile rank score test, which was proposed by Gutenbrunner et al. [6] for linear models.

We partition $\beta$ into two parts $\beta_1 \in \mathbb{R}^{q_1}$ and $\beta_2 \in \mathbb{R}^{q_2}$ with $q_1 + q_2 = q$. Suppose we want to test $H_0 : \beta_1 = 0$. Let $Z^{(1)}$ be the $N \times q_1$ and $Z^{(2)}$ be the $N \times q_2$ design matrices corresponding to $\beta_1$ and $\beta_2$, respectively. Furthermore, we denote $\mathcal{W} = (\Pi, Z^{(2)})$, $P_\varpi = \mathcal{W}(\mathcal{W}'B\mathcal{W})^{-1}\mathcal{W}B$, $\mathcal{D} = (I - P_\varpi)Z^{(1)}$, and $\phi = (\Theta', \beta_2')'$. In practice, the weight matrix $B$ can be estimated by $\hat{B}$ as defined in (3.2), and this will not affect the asymptotic behavior of the rank score test statistics to be proposed in the following and in Section 4.

Let $d_{ij} \in \mathbb{R}^{q_1}$ and $\varpi_{ij} \in \mathbb{R}^{p_{k_n}+q_2}$ be the column components of $\mathcal{D}$ and $\mathcal{W}$ associated with the $j$th measurement of the $i$th subject, respectively, and $\mathcal{D}_i = (d_{i1}, \ldots, d_{im_i})'$. We define the rank score test statistic as

$$\mathcal{T}_n = S_n' \hat{V}_n^{-1} S_n, \tag{3.4}$$

where

$$S_n = N^{-1/2} \sum_{ij} d_{ij} \psi_\tau(\hat{e}_{ij}), \qquad \hat{e}_{ij} = y_{ij} - \varpi_{ij}'\hat{\phi},$$

$$\hat{\phi} = \arg\min_{\phi \in \mathbb{R}^{p_{k_n}+q_2}} \sum_{ij} \rho_\tau(y_{ij} - \varpi_{ij}'\phi), \qquad \hat{V}_n = N^{-1} \sum_{i=1}^n \mathcal{D}_i' \psi_\tau(\hat{e}_i) \psi_\tau'(\hat{e}_i) \mathcal{D}_i$$

and

$$\psi_\tau(\hat{e}_i) = (\psi_\tau(\hat{e}_{i1}), \ldots, \psi_\tau(\hat{e}_{im_i}))'.$$

To establish the asymptotic distribution of $\mathcal{T}_n$, we modify the assumption A1 as A1' and make an additional assumption A7,

A1'  There exists some $r > 2$ such that $\alpha_l(t) \in \mathcal{H}_r, l = 1, \ldots, p$.
A7  The minimum eigenvalue of $V_n \doteq N^{-1} \sum_{i=1}^n \mathcal{D}_i' A_i(\Delta) \mathcal{D}_i$ is bounded away from zero for sufficiently large $n$.

THEOREM 3.  *If assumptions* A1' *and* A2–A7 *hold, and* $n^{1/(4r)} \ll k_n \ll n^{1/4}$ *with* $a_n \ll b_n$ *meaning* $a_n = o(b_n)$, *we have* $\mathcal{T}_n \xrightarrow{D} \chi^2(q_1)$ *as* $n \to \infty$.

REMARK 1.  The finite sample efficiency of $\hat{V}_n$ could be debilitated by imposing a fully nonparametric structure to the correlation matrix. If the structure of $A_i(\Delta)$ is known, we can estimate $\Delta$ empirically with $\hat{\Delta}$ by incorporating information across subjects. Plugging $\hat{\Delta}$ to $V_n$, we obtain an asymptotically equivalent estimator $\hat{V}_{2n} = N^{-1} \sum_{i=1}^n \mathcal{D}_i^T A(\hat{\Delta}) \mathcal{D}_i$. Suppose, for instance, that $A_i(\Delta)$ has a compound symmetry structure with $\delta_{ij_1j_2} = \delta$. The parameter $\delta$ can be estimated consistently by $\hat{\delta} = L^{-1} \sum_i \sum_{j_1 \neq j_2} I(\hat{e}_{ij_1} < 0, \hat{e}_{ij_2} < 0)$, where $L$ denotes the total number of pairs of repeated measurements from the same subject. For this special case, as $\hat{V}_{2n}$ involves the estimation of only one nuisance parameter $\delta$, it is expected to be more efficient than $\hat{V}_n$ in finite samples. The same discussion applies to the estimation of the covariance matrix of $s_n$ in (4.4).

**4. Constancy of varying coefficients.** In semiparametric models, another question of practical interest is to test whether one or some of the varying coefficients is constant. We propose a rank-score-based procedure for testing such hypotheses through model re-parameterization. Without loss of generality, consider testing the null hypothesis that the first $1 \le p_1 \le p$ coefficient functions $\alpha_l(\cdot)$ are constant:

$$H_0 : \alpha_l(t) = \gamma_l, \qquad l = 1, \ldots, p_1, \tag{4.1}$$

for all $t \in [0,1]$, where $\gamma_l$ are unknown constants, versus the alternative hypothesis

(4.2) $\quad H_1$ : one or more functions $\alpha_l(t), l = 1, \ldots, p_1,$ are time-varying.

We can find transformation matrices $G$ and $\widetilde{G}$ such that

$$G\pi_{ij} = (1, \bar{\pi}_{ij}')', \qquad \widetilde{G}\Pi_{ij} = (\Pi_{ij}^{(1)}, \Pi_{ij}^{(2)}),$$

where $\Pi_{ij}^{(1)} = (x_{ij}^{(1)}\bar{\pi}_{ij}', \ldots, x_{ij}^{(p_1)}\bar{\pi}_{ij}')' \in \mathbb{R}^{\tilde{p}_1}$ and $\Pi_{ij}^{(2)} \in \mathbb{R}^{p k_n - \tilde{p}_1}$, $\tilde{p}_1 = p_1(k_n + \hbar)$. Let

$$\Pi^{(1)} = (\Pi_{11}^{(1)}, \ldots, \Pi_{1m_1}^{(1)}, \ldots, \Pi_{nm_n}^{(1)})', \qquad \Pi^{(2)} = (\Pi_{11}^{(2)}, \ldots, \Pi_{1m_1}^{(2)}, \ldots, \Pi_{nm_n}^{(2)})'.$$

Denote $\xi_1 \in \mathbb{R}^{\tilde{p}_1}$ and $\xi_2 \in \mathbb{R}^{p k_n - \tilde{p}_1}$ as the parameters corresponding to the design matrices $\Pi^{(1)}$ and $\Pi^{(2)}$, respectively. The $\tau$th conditional quantile of $y_{ij}$ can then be approximated by

$$Q_{y_{ij}}(\tau | x_{ij}, z_{ij}, t_{ij}) = \Pi_{ij}^{(1)\prime}\xi_1 + \Pi_{ij}^{(2)\prime}\xi_2 + z_{ij}'\beta \tag{4.3}$$

and (4.1) and (4.2) can be represented as $\widetilde{H}_0 : \xi_1 = 0$ versus $\widetilde{H}_1 : \xi_1 \ne 0$.

The rank score test is based on the quantile estimates $\hat{\varphi} = (\hat{\xi}_2', \hat{\beta}')'$ of $\varphi = (\xi_2', \beta')'$ obtained under $H_0$. Let $W = (\Pi^{(2)}, Z)$, $P_w = W(W'BW)^{-1}W'B$, $D = (I - P_w)\Pi^{(1)}$, and let $w_{ij}$ and $d_{ij}$ be the column components of $W$ and $D$ associated with the $j$th measurement of the $i$th subject, respectively. Furthermore, we denote $D_i = (d_{i1}, \ldots, d_{im_i})'$. Note that the dimensions of $w_{ij}$ and $d_{ij}$ are both increasing in the dimension of the B-spline space. We denote the rank score test statistic as

$$t_n = s_n' \hat{v}_n^{-1} s_n, \tag{4.4}$$

where $s_n = N^{-1/2} \sum_{i=1}^n \sum_{j=1}^{m_i} d_{ij} \psi_\tau(\hat{e}_{ij})$, $\hat{e}_{ij} = y_{ij} - \Pi_{ij}^{(2)\prime}\hat{\xi}_2 - z_{ij}'\hat{\beta} = y_{ij} - w_{ij}'\hat{\varphi}$, $\hat{v}_n = N^{-1} \sum_{i=1}^n D_i'\psi_\tau(\hat{e}_i)\psi_\tau(\hat{e}_i)'D_i$ and $\psi_\tau(\hat{e}_i) = (\psi_\tau(\hat{e}_{i1}), \ldots, \psi_\tau(\hat{e}_{im_i}))'$.

If $k_n = k$ is bounded corresponding to the parametric model (4.3), a rank score test based on $t_n$ and the $\chi^2(p_1(k + \hbar))$ reference distribution can be used to test $\widetilde{H}_0 : \xi_1 = 0$; see [24]. However, to test the constancy of functional coefficients in the PLVC models, this rank score test will be inconsistent unless $k_n$ grows with the sample size. The following Theorem 4 gives the asymptotic null distribution of the rank score test statistic $t_n$ for growing $k_n$.

THEOREM 4. *Assume that* A1$'$ *and* A2–A6 *hold*, $\hbar \geq 3$, *and the number of knots satisfies* $n^{1/(2r+2)} \ll k_n \ll n^{1/5}$. *Furthermore, assume that the minimum eigenvalue of the matrix* $v_n \doteq N^{-1}\sum_{i=1}^{n} D_i' A_i(\Delta) D_i$ *is bounded away from zero for sufficiently large* $n$. *Under* $H_0$, *we have*

$$\frac{t_n - (k_n + \hbar)p_1}{\sqrt{2(k_n + \hbar)p_1}} \xrightarrow{D} N(0,1) \qquad as\ k_n \to \infty.$$

REMARK 2. In the special homoscedastic case, where the errors have a common density with $f_{ij}(0) = f(0)$ for all $i$ and $j$, the weight matrix $B$ will be canceled out in the projection matrices such as $P_\varpi$ and $P_w$. Therefore, the rank score test statistics (3.4) and (4.4) reduce to simpler forms free of $f$. However, for heteroscedastic errors, the matrix $B$ needs to be incorporated appropriately in the projection matrices to ensure the orthogonality of $\mathcal{D}$ and $B\mathcal{W}$, $D$ and $BW$, and the sandwich formula in the covariance matrix of $\hat{\beta}$.

Kim [15] proposed a similar inference procedure to test whether all of the coefficients are constant. Our developed testing procedure has wider applications, as it can examine the constancy of a subset of possibly varying coefficients, and thus can be used for model selection. Note that in Kim [15], the dimension of nuisance parameters under $H_0$, $\xi_2$, is finite. In our setup, both the dimension of nuisance parameters and that of parameters to be tested, $\xi_1$, are allowed to increase with sample size $n$. The increasing dimension of parameters poses challenges to establish the asymptotic representation of $s_n$, which is needed to prove Theorem 4.

**5. Simulation study.** In this section, we investigate the finite sample performance of the proposed estimation and inference methods with Monte Carlo simulation studies. We generate 2000 data sets, each consisting of $n$ subjects. We consider two sample sizes $n = 30$ and $n = 100$. Each subject is supposed to be measured at scheduled time points $\{0, 1, 2, \ldots, 10\}$, each of which (excluding time 0) has a 20% probability of being skipped. This results in different numbers of repeated measurements $m_i$ for each subject. The actual measurement times are generated by adding an $U(-0.5, 0.5)$ random variable to the nonskipped scheduled times. We focus on $\tau = 0.25$ and $\tau = 0.5$ in this study.

The data sets are generated from the following heteroscedastic model

$$\begin{aligned}(5.1)\qquad y_{ij} &= \alpha_0(t_{ij}) + \alpha_1(t_{ij})x_{ij,1} + \alpha_2(t_{ij})x_{ij,2} + \alpha_3(t_{ij})x_{ij,3} \\ &\quad + \beta z_i + (1 + |x_{ij,1}|)e_{ij}(\tau),\end{aligned}$$

where $x_{ij,1}$, $x_{ij,2}$ and $x_{ij,3}$ follow the standard normal, $U(t_{ij}/10, 2 + t_{ij}/10)$, and the standard exponential distributions, respectively, $z_i$ is a time-invariant covariate generated from a Bernoulli distribution with 0.5 probability of being 1, $e_{ij}(\tau) = e_{ij} - F^{-1}(\tau)$ with $F$ being the common CDF of $e_{ij}$. Here, $F^{-1}(\tau)$ is subtracted from $e_{ij}$ to make the $\tau$th quantile of $e_{ij}(\tau)$ zero for identifiability purpose.

We consider three cases for generating $e_{ij}$. In cases 1 and 2, $e_i = (e_{i1}, \ldots, e_{im_i})$ follows a multivariate normal distribution $N(0, \Sigma_i)$. In case 3, $e_i$ is generated from a multivariate $t(3)$ distribution, specifically $e_i = \sqrt{3} u^{-1/2} \varepsilon$, where $\varepsilon \sim N(0, \Sigma_i)$, $u$ consists of $m_i$ independent $\chi^2(3)$ random variables, and $\varepsilon$ and $u$ are mutually independent. The covariance matrix $\Sigma_i$ has an exchangeable correlation structure with correlation coefficient of $\varrho = 0.8$ in cases 1 and 3, while an AR(1) correlation structure with $\varrho(e_{ij_1}, e_{ij_2}) = 0.8^{|t_{ij_1} - t_{ij_2}|}$ in case 2. We set $\beta = 1$ and

$$\alpha_0(t) = 15 + 20 \sin\left(\frac{t\pi}{20}\right),$$

(5.2) $$\alpha_1(t) = 2 - 3\cos\left(\frac{(3t-25)\pi}{15}\right),$$

$$\alpha_2(t) = 6 - 0.6t, \qquad \alpha_3(t) = -4 + (20 - 3t)^3/1000.$$

5.1. *Estimation.* Throughout our numerical studies, we use the cubic splines with $\hbar = 3$ in the B-spline approximation. For each data set, we choose $k_n$ as the minimizer to the following Schwarz-type information criterion,

(5.3)
$$\text{SIC}(k) = \log\left\{\sum_{ij} \rho_\tau(y_{ij} - \Pi'_{ij}\hat{\theta}_{(k)} - z'_i\hat{\beta}_{(k)})\right\} + \frac{\log N}{2N}\{p(\hbar + k + 1) + q\},$$

where $p = 4$, $q = 1$, $\hat{\theta}_{(k)}$ and $\hat{\beta}_{(k)}$ are the $\tau$th quantile estimators obtained from minimizing (2.1) with $k$ knots; see He, Zhu and Fung [10] for a similar criterion for knots selection.

To demonstrate the flexibility and efficiency of the PLVC model, we compare its performance to the linear constant coefficient (LCC) model

$$y_{ij} = \alpha_0 + \alpha_1 x_{ij,1} + \alpha_2 x_{ij,2} + \alpha_3 x_{ij,3} + \beta z_i + e_{ij}.$$

Let $\hat{\beta}_{\text{PLVC}}$ and $\hat{\beta}_{\text{LCC}}$ denote the quantile estimation of $\beta$ obtained under the PLVC and the LCC models, respectively. For fair comparison, we also generate 2000 data sets from the LCC model with $(\alpha_0, \alpha_1, \alpha_2, \alpha_3) = (15, 2, 6, -4)$ for each case. Table 1 summarizes the estimated mean squared error (MSE) and bias of $\hat{\beta}_{\text{PLVC}}$ and $\hat{\beta}_{\text{LCC}}$ at $\tau = 0.5$ and $n = 100$. When the data is generated from a PLVC model, fitting the LCC model leads to less efficient estimations. The PLVC model is more flexible and it does not lose much efficiency even when the underlying model has constant coefficients.

5.2. *Inference on $\beta$.* We vary $\beta$ in the PLVC model (5.1) from 0 to 2.5 to assess the type I error and power of the proposed tests. We consider two variants of the proposed rank score test: QRS and $\text{QRS}_\delta$. The QRS does not specify the dependence structure of errors and is based on the covariance estimator $\hat{V}_n$. The $\text{QRS}_\delta$ is based on $\hat{V}_{2n}$ assuming an exchangeable intra-subject correlation. For comparison, we also include the Wald-type test (Wald), and the mean regression method of Huang, Wu and Zhou [13] based on bootstrap of 500 samples, referred to as LSE. The comparison between LSE and the other three methods is possible only at $\tau = 0.5$ when the mean and median functionals are the same.

Table 2 summarizes the type I errors of the four methods in cases 1–3. The Wald-type test relies heavily on the estimation of the error density function. For all the cases considered, Wald fails to maintain the nominal significance level of 0.05, especially at $n = 30$. The other three tests all maintain the level reasonably well. Figure 1 shows the power curves of QRS, $\text{QRS}_\delta$ and LSE in cases 1 and 3. The rank score tests QRS and $\text{QRS}_\delta$ give very similar power at both $n = 30$ and $n = 100$ in this simulation study for inference on the scalar $\beta$. Another power comparison (not reported due to space limit) suggests that $\text{QRS}_\delta$ is more powerful than QRS in small samples such as $n = 30$ for inference on $\beta \in \mathbb{R}^q$ with $q > 1$. Compared to LSE, the median approaches are much more powerful for heavy-tailed data (case 3), and they perform equally well in finite samples for normal data (cases 1–2), where LSE is known to have higher asymptotic efficiency.

TABLE 1
*The Monte Carlo mean squared error (MSE) and bias of $\hat{\beta}_{\text{PLVC}}$ and $\hat{\beta}_{\text{LCC}}$ in cases 1–3 at $\tau = 0.5$ and $n = 100$*

|  | Underlying model: PLVC | | | | Underlying model: LCC | | | |
|---|---|---|---|---|---|---|---|---|
|  | MSE | | Bias | | MSE | | Bias | |
|  | $\hat{\beta}_{\text{PLVC}}$ | $\hat{\beta}_{\text{LCC}}$ | $\hat{\beta}_{\text{PLVC}}$ | $\hat{\beta}_{\text{LCC}}$ | $\hat{\beta}_{\text{PLVC}}$ | $\hat{\beta}_{\text{LCC}}$ | $\hat{\beta}_{\text{PLVC}}$ | $\hat{\beta}_{\text{LCC}}$ |
| Case 1 | 0.111 | 0.187 | 0.001 | $-0.001$ | 0.111 | 0.111 | 0.001 | 0.001 |
| Case 2 | 0.073 | 0.160 | 0.001 | $-0.007$ | 0.072 | 0.073 | $-0.002$ | $-0.003$ |
| Case 3 | 0.156 | 0.336 | $-0.014$ | $-0.016$ | 0.154 | 0.154 | $-0.016$ | $-0.014$ |

TABLE 2
*Type I errors for testing $H_0 : \beta = 0$. The nominal significance level is 0.05*

| Case | $\tau$ | $n = 30$ | | | | $n = 100$ | | | |
|---|---|---|---|---|---|---|---|---|---|
|  |  | QRS | $\text{QRS}_\delta$ | Wald | LSE | QRS | $\text{QRS}_\delta$ | Wald | LSE |
| 1 | 0.25 | 0.059 | 0.065 | 0.139 | / | 0.061 | 0.061 | 0.113 | / |
|  | 0.5 | 0.061 | 0.068 | 0.130 | 0.069 | 0.053 | 0.053 | 0.088 | 0.053 |
| 2 | 0.25 | 0.059 | 0.062 | 0.145 | / | 0.055 | 0.057 | 0.117 | / |
|  | 0.5 | 0.057 | 0.068 | 0.104 | 0.072 | 0.059 | 0.059 | 0.097 | 0.056 |
| 3 | 0.25 | 0.064 | 0.063 | 0.129 | / | 0.049 | 0.049 | 0.075 | / |
|  | 0.5 | 0.054 | 0.064 | 0.100 | 0.038 | 0.049 | 0.050 | 0.081 | 0.030 |

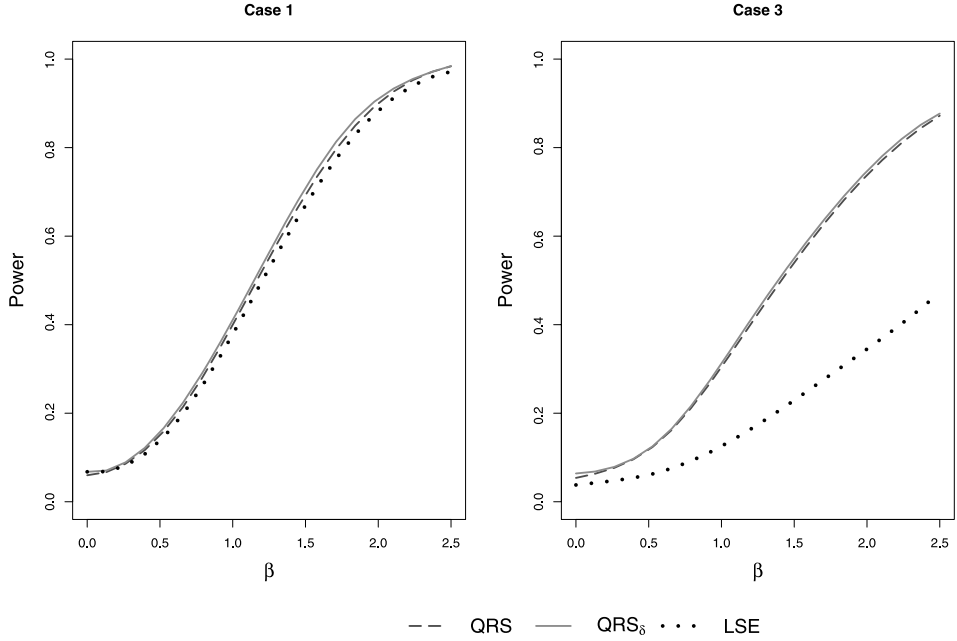

Fig. 1. *Power curves of different methods for testing $H_0: \beta = 0$ at $\tau = 0.5$ and $n = 30$.*

5.3. *Inference on the constancy of $\alpha_1(t)$.* To test whether $\alpha_1(t)$ deviates from a constant, we fix $\beta = 1$ and generate data from a sequence of alternative models

$$\alpha_1(t, \eta) = 2 - 3\eta \cos\left(\frac{(t - 25)\pi}{15}\right),$$

where $\eta$ determines the extent that $\alpha_1(t)$ varies with time. We set $\eta \in [0, 1.5]$ with 0 corresponding to a model with a constant coefficient $\alpha_1$.

Table 3 shows that QRS, $\text{QRS}_\delta$ and LSE all maintain the level reasonably well, even at small samples with $n = 30$. The QRS and $\text{QRS}_\delta$ give similar power at $n = 100$, but $\text{QRS}_\delta$ is more powerful in small samples $n = 30$; see Figure 2 for typical examples in cases 1 and 3 at $\tau = 0.5$. As observed in Section 5.2, the median approaches perform competitively with LSE in cases 1–2, and they are substantially more powerful than LSE in case 3.

We now summarize the main findings from this simulation study. First, the proposed rank score tests perform well in terms of both level and power, for testing $\beta$ or the constancy of varying coefficients. Second, for small samples, $\mathrm{QRS}_\delta$ is more powerful than QRS for testing the constancy of $\alpha_1(t)$ and for testing $\beta \in \mathbb{R}^q$ with $q > 1$ in cases 1 and 3 when the underlying correlation is exchangeable. For the data sets generated in this study, $\mathrm{QRS}_\delta$ performs competitively even in case 2 when the correlation structure is misspecified. It is commonly known that analysis at different quantiles can reveal structural heterogeneity that might be overlooked by mean regression methods. When the main interest is on the center of the distribution, the proposed median approach does not lose much finite sample efficiency for normal data compared to the mean method, but it performs more robustly for heavy tailed data.

## 6. Application to AIDS data.

6.1. *Background of the study.* We illustrate the proposed quantile approach by analyzing a subset from the Multi-Center AIDS Cohort study. The data set consists of 283 homosexual males who were infected with HIV between 1984 and 1991. Each patient had a different number of repeated measurements and different measuring times. Details of the experimental design can be found in Kaslow et al. [14]. Several researchers including Fan, Huang and Li [5], Huang, Wu and Zhou [13] and Qu and Li [21] have studied the same data set to analyze the mean CD4 percentage. Our analysis aims to model the effects of smoking status, age and pre-HIV infection CD4 percentage (PreCD4) on different quantiles of the CD4 percentage population.

Previous analyses on the same data set suggested that smoking status and age have no significant effects on the mean CD4 percentage, and the baseline has time-varying effect, but it remains unclear whether the effect of PreCD4 is constant or not. Therefore, at a given quantile level $\tau$, we consider a parsimonious PLVC model

$$(6.1) \quad y_{ij} = \alpha_0(t_{ij}, \tau) + \alpha_1(t_{ij}, \tau)x_i + \beta_1(\tau)z_{i,1} + \beta_2(\tau)z_{i,2} + e_{ij}(\tau),$$

Table 3
*Type I errors for testing the constancy of $\alpha_1(t)$. The nominal significance level is 0.05*

| Case | $\tau$ | $n = 30$ | | | $n = 100$ | | |
|---|---|---|---|---|---|---|---|
| | | QRS | $\mathrm{QRS}_\delta$ | LSE | QRS | $\mathrm{QRS}_\delta$ | LSE |
| 1 | 0.25 | 0.037 | 0.061 | / | 0.060 | 0.059 | / |
|   | 0.5  | 0.052 | 0.057 | 0.064 | 0.049 | 0.048 | 0.056 |
| 2 | 0.25 | 0.035 | 0.053 | / | 0.055 | 0.059 | / |
|   | 0.5  | 0.057 | 0.062 | 0.070 | 0.051 | 0.051 | 0.067 |
| 3 | 0.25 | 0.036 | 0.055 | / | 0.047 | 0.046 | / |
|   | 0.5  | 0.044 | 0.056 | 0.021 | 0.061 | 0.065 | 0.022 |

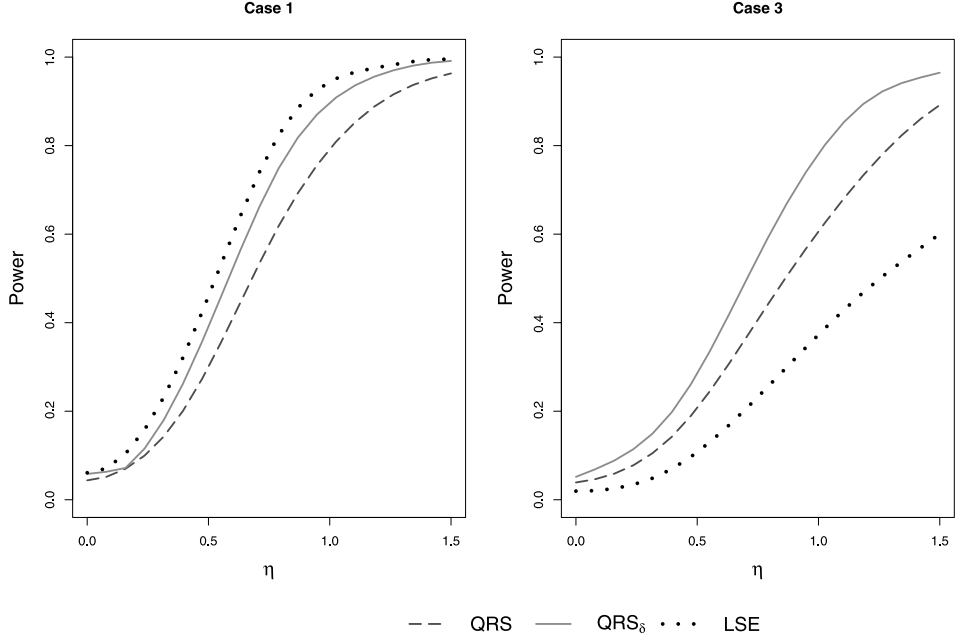

Fig. 2. *Power curves of different methods for testing the constancy of $\alpha_1(t)$ at $\tau = 0.5$ and $n = 30$.*

where $y_{ij}$ is the $i$th patient's CD4 cell percentage at time $t_{ij}$ (in years), $x_i$ is the centered PreCD4, $z_{i,1}$ is the $i$th patient's smoking status (1 for smoker and 0 for nonsmoker), and $z_{i,2}$ is the $i$th patient's centered age at HIV infection. In this study, we consider a set of quantiles with $\tau \in \{0.1, 0.2, \ldots, 0.8, 0.9\}$. At a fixed quantile level $\tau$, the baseline function $\alpha_0(t, \tau)$ represents the $\tau$th quantile of CD4 percentage $t$ years after the infection for a nonsmoker with average PreCD4 percentage and average age at HIV infection.

6.2. *Model assessment.* For comparative purpose, we also consider the LCC model assuming that both $\alpha_0(t, \tau)$ and $\alpha_1(t, \tau)$ are constant at a given quantile level $\tau$. To assess how well the two models fit the data, we consider the following model assessment tool by comparing the empirical distribution of $Y$ with the simulated distribution from each model. We first generate $u$ from $U(0, 1)$. At a fixed time point $t^*$, we randomly choose a subject from those who have measurements at $t^*$ or close to $t^*$ (within a distance of 0.001), and denote the corresponding covariates as $(x^*, z_1^*, z_2^*)$. The simulated $Y^*$ is generated as the $u$th conditional quantile $\hat{\alpha}_0(t^*, u) + \hat{\alpha}_1(t^*, u)x^* + \hat{\beta}_1(u)z_1^* + \hat{\beta}_2(u)z_2^*$, where $(\hat{\alpha}_0(\cdot, u), \hat{\alpha}_1(\cdot, u), \hat{\beta}_1(u), \hat{\beta}_2(u))$ are the estimated $u$th quantile coefficients in the model under assessment. Repeating this procedure many times, say 500, we can obtain a simulated sample. If the model fits data well, the marginal distribution of the simulated $Y^*$ should match that of the observed $Y$. Figure 3 shows the Q–Q plots of the empirical $(Y)$ and simulated CD4 percentage $(Y^*)$ at $t^* = 0.2, 0.8, 1.2$ and 3.8. The Q–Q plots suggest that the PLVC model fits the data universally better than the constant coefficient model.

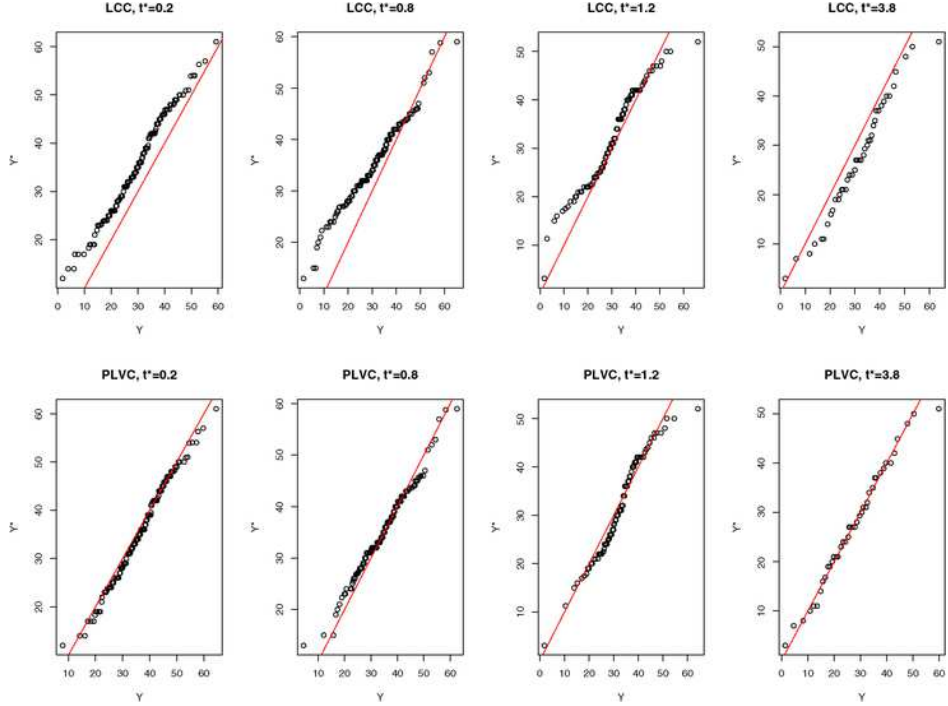

Fig. 3. *Assessment of the PLVC and LCC models for fitting the CD4 data. The diagonal line in each plot is $y = x$.*

6.3. *Hypothesis testing.* At each quantile level $\tau$, we consider four null hypotheses, $H_{01}$: the baseline curve $\alpha_0(t, \tau)$ is time invariant; $H_{02}$: the PreCD4 effect is time invariant; $H_{03}$: smoking has no effect and $H_{04}$: age has no effect. We apply the proposed quantile rank score test $\text{QRS}_\delta$ by assuming an exchangeable correlation structure to test each of the four null hypotheses. The $p$-values are summarized in the first four rows of Table 4. Figure 4(a) and (b) plot the estimated baseline and PreCD4 effects at $\tau = 0.1, 0.3$ and 0.7.

Consistent with the mean regression results, our analysis indicates that neither smoking nor age has any significant effects, and the baseline curve is significantly time-varying at all quantile levels considered. Additionally, Figure 4(a) suggests that the CD4 percentage in the lower quantiles, that is, for the group of baseline people with more severe illness, depletes rather quickly across the time that was considered. In contrast, the upper quantiles of the baseline CD4 percentage drop at a slower rate for the first 3 years, and become stable afterward.

There is some disagreement in terms of hypothesis testing on the constancy of the PreCD4 effect in literature. Huang, Wu and Zhou [13] reported a $p$-value of 0.059 for testing $H_{02}$, while Qu and Li [21] claimed that the effect of PreCD4 is significantly time-varying with a $p$-value of 0.045. Our proposed quantile regression method suggests that the effect of PreCD4 is time-decaying with marginal significance at lower quantiles $\tau = 0.2$ and $0.3$, and it tends to be stable for $\tau > 0.3$. Such structural heterogeneity would not have been revealed by ordinary regression methods focusing solely on the conditional mean.

To examine the constancy of functional coefficients, we also investigate an alternative shrinkage approach through minimizing the penalized objective function: $\sum_{ij} \rho_\tau(y_{ij} - \Pi_{ij}^{(1)\prime}\xi_1 - \Pi_{ij}^{(2)\prime}\xi_2 - z'_{ij}\beta) + \lambda\|\xi_1\|_1$, where $\|\cdot\|_1$ denotes the $L_1$ norm, and $\xi_1$ and $\xi_2$ are defined in (4.3). The tuning parameter $\lambda$ is chosen by minimizing the Schwarz-type information criterion SIC = $\log\{\sum_{ij}\rho_\tau(y_{ij} - \Pi_{ij}^{(1)\prime}\hat{\xi}_1 - \Pi_{ij}^{(2)\prime}\hat{\xi}_2 - z'_{ij}\hat{\beta})\} + \frac{\log N}{2N}df$, where $df$ is the number of coefficients that are not shrunk to zero. If all components of $\xi_1$ are shrunk to zero simultaneously, that is, $\|\xi_1\|_1 = 0$, the coefficient functions $\alpha_l(t), l = 1, \ldots, p_1$, are considered as time-invariant. We apply the shrinkage approach to examine the constancy of the baseline and the PreCD4 effects separately. The resulting $\|\hat{\xi}_1\|_1$ at different quantile levels are summarized in the last two rows of Table 4. The results from the shrinkage method agree quite well

TABLE 4
*Analysis results for the CD4 data at different quantile levels. The first four rows are the $p$-values from the rank score tests. The last two rows are the resulting $\|\hat{\xi}_1\|_1$ from the shrinkage approach for examining the constancy of the baseline and PreCD4 effects, respectively*

| Null hypothesis | $\tau$ | | | | | | |
|---|---|---|---|---|---|---|---|
| | 0.1 | 0.2 | 0.3 | 0.4 | 0.5 | 0.7 | 0.9 |
| | The $p$-values from the rank score tests | | | | | | |
| $H_{01}$: constant baseline | 0.000 | 0.000 | 0.000 | 0.000 | 0.000 | 0.000 | 0.000 |
| $H_{02}$: constant PreCD4 | 0.264 | 0.028 | 0.010 | 0.332 | 0.301 | 0.324 | 0.597 |
| $H_{03}$: no smoking effect | 0.085 | 0.054 | 0.476 | 0.685 | 0.807 | 0.581 | 0.211 |
| $H_{04}$: no age effect | 0.631 | 0.525 | 0.214 | 0.117 | 0.212 | 0.268 | 0.418 |
| | The $\|\hat{\xi}_1\|_1$ from the shrinkage approach | | | | | | |
| $H_{01}$: constant baseline | 0.781 | 0.571 | 0.864 | 0.663 | 0.484 | 0.687 | 0.623 |
| $H_{02}$: constant PreCD4 | 0.295 | 0.284 | 0.153 | 0.000 | 0.345 | 0.000 | 0.000 |

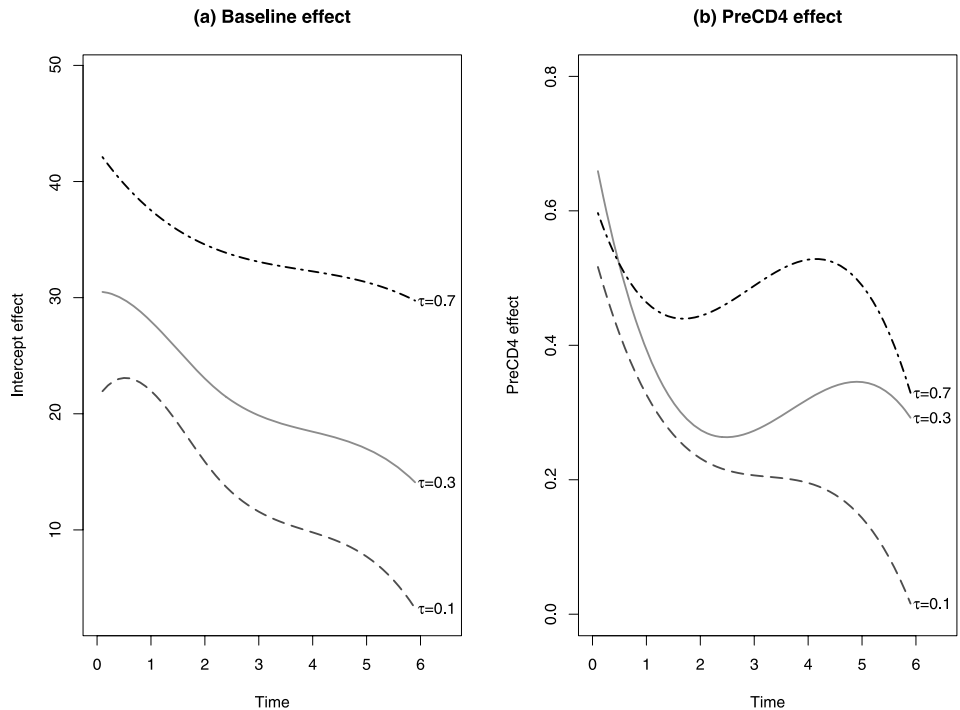

Fig. 4. *The estimated varying coefficient curves at $\tau = 0.1, 0.3$ and $0.7$.*

with those from the rank score test, both suggesting that the baseline effect is time-varying at all quantiles considered, and the PreCD4 effect tends to be time-varying at $\tau = 0.2$ and $0.3$. Further investigation is needed to understand the theoretical property of the shrinkage method and to compare with the proposed rank score test approach.

**7. Discussion.** In this paper, we introduced a marginal quantile partially linear model with varying coefficients for analyzing longitudinal data. We proposed a simple procedure for estimating the quantile coefficients $\beta$ and $\alpha(t)$ using B-spline approximation, and established the asymptotic properties of the resulting estimators. We further developed rank score tests for some important inferential issues on both $\beta$ and $\alpha(t)$. The proposed rank score tests are easy to implement and robust in performance. The estimating and inference approach does not require any specification of the error distribution or the dependence structure. However, our empirical studies suggest that the finite sample performance of the tests could be improved by specifying an approximate correlation structure.

We chose B-spline to approximate the smooth functional coefficients in this paper due to its computational efficiency and stability. The B-spline is able to exhibit local features and it provides good approximations with a small number of knots. Other smoothing techniques, such as penalized regression splines and local polynomial fitting could be employed in our model framework, and this deserves further research and comparison.

The proposed rank score tests can be inverted to construct confidence intervals for the parametric coefficient $\beta$; see Koenker [16] for details of such confidence interval construction in a different context. However, we have not addressed the issue of confidence regions for the functional curves $\alpha_l(t)$, such as intervals at a fixed $t$ or simultaneous confidence bands for all $t$ within a range. Resampling of subjects might provide a solution to this problem, but further theoretical and practical investigation is needed.

## APPENDIX

Throughout the appendix, we use $\|\cdot\|$ to denote the Euclidean norm. First note that by assumption A1 and Corollary 6.21 of [22], there exists a spline approximation $\pi'(t)\theta_{l0}$ to $\alpha_l(t)$ such that

$$(A.1) \qquad \sup_{t \in [0,1]} |\alpha_l(t) - \pi'(t)\theta_{l0}| = O(k_n^{-r}), \qquad l = 1, \ldots, p.$$

To build links between $(\hat{K}_n, \hat{\Lambda}_n)$ and $(K_n, \Lambda_n)$, we define

$$P = \Pi(\Pi'B\Pi)^{-1}\Pi'B, \qquad Z^* = (I - P)Z,$$

$$K_n^* = Z^{*\prime}BZ^* = \sum_{i=1}^n Z_i^{*\prime}B_iZ_i^*, \qquad \Lambda_n^* = Z^{*\prime}AZ^* = \sum_{i=1}^n Z_i^{*\prime}A_iZ_i^*.$$

LEMMA 1. *Under the assumptions of Theorem 1, we have*

$$N^{-1}(K_n^* - K_n) = o_p(1) \quad \text{and} \quad N^{-1}(\Lambda_n^* - \Lambda_n) = o_p(1).$$

PROOF. Let $\Phi^* = (\phi^*(X_1, T_1)', \ldots, \phi^*(X_n, T_n)')'$ and $Z = (Z - \Phi^*) + \Phi^* \doteq \Delta + \Phi^*$, where $\phi^*$ is defined in (2.2). We can write

$$(A.2) \qquad \begin{aligned} N^{-1}K_n^* = N^{-1}\{&\Delta'B\Delta + \Phi^{*\prime}(I - P')B(I - P)\Phi^* \\ &+ \Delta'(I - P')B(I - P)\Phi^* \\ &+ \Phi^{*\prime}(I - P')B(I - P)\Delta - \Delta'P'BP\Delta\}. \end{aligned}$$

By (A.1) and assumption A5, there exists a matrix $M$ such that $\|\Phi^* - \Pi M\| = O_p(n^{1/2}k_n^{-r})$. In addition, as $P$ is a projection matrix, we have

$$\|(I - P)\Phi^*\| = \|\Phi^* - \Pi M\| + \|\Pi M - P\Phi^*\| \leq 2\|\Phi^* - \Pi M\| = O_p(n^{1/2}k_n^{-r}).$$

Let $m(X,T) = (m(X_1,T_1)',\ldots,m(X_n,T_n)')'$, $\eta = Z - m(X,T)$ and $\zeta = m(X,T) - \Phi^*$. By assumption A5, we have

(A.3) $\quad \|P\eta'\|^2 = O_p(\mathrm{tr}(P'P)) = O_p(k_n), \qquad \|(I-P)\eta\| = O_p(n^{1/2}).$

On the other hand, since $\phi_l^*(\cdot,\cdot)$ is the projection of $m_l(\cdot,\cdot)$ onto the varying coefficient functional space $\mathcal{Y}$ (i.e., $\zeta \perp \mathcal{Y}$) and $\Pi \in \mathcal{Y}$, we have

$$E(\Pi'B\zeta) = 0 \quad \text{and} \quad E\|\Pi'B\zeta\|^2 = E\left(\sum_{i=1}^n \|\Pi_i B_i \zeta_i\|^2\right) = O_p(k_n n),$$

where $\Pi_i$ and $\zeta_i$ are the rows of $\Pi$ and $\zeta$, respectively, corresponding to the $i$th subject. Hence $\|\Pi'B\zeta\| = O_p(k_n^{1/2} n^{1/2})$. By Lemma A.4 of [15], we have

$$\|P\zeta\| \le \|\Pi\| \|(\Pi'B\Pi)^{-1}\| \|\Pi B\zeta\|$$
$$= O_p(n^{1/2} k_n^{-1/2}) O_p(k_n/n) O_p((k_n n)^{1/2}) = O_p(k_n),$$

which together with (A.3) gives

(A.4) $\qquad\qquad\qquad \|P\Delta\| = O_p(k_n).$

Thus, all the last four terms on the right-hand side of (A.2) are $o_p(1)$, so the result is proven as desired. The proof of $N^{-1}(\Lambda_n^* - \Lambda_n) = o_p(1)$ follows with similar arguments. $\square$

LEMMA 2. *Under the assumptions of Theorem 1, we have*

(A.5) $\qquad\qquad \Lambda_n^{-1/2} Z^{*\prime} \psi_\tau(e) \longrightarrow N(0, I_q),$

*where $\psi_\tau(e) = (\psi_\tau(e_{11}),\ldots,\psi_\tau(e_{nm_n}))'$.*

PROOF. Similar to (A.2), we can write

(A.6) $\quad \begin{aligned} \Lambda_n^{-1/2} Z^{*\prime} \psi_\tau(e) &= \Lambda_n^{-1/2} \Delta' \psi_\tau(e) \\ &\quad + \Lambda_n^{-1/2}\{\Phi^{*\prime}(I-P')\psi_\tau(e) - \Delta'P'\psi_\tau(e)\}. \end{aligned}$

By (A.1) and (A.4), we have $N^{-1} E\|\Phi^{*\prime}(I-P')\psi_\tau(e)\|^2 = o(1)$, and $N^{-1} \times E\|\Delta'P'\psi_\tau(e)\|^2 = o(1)$. Thus, the last term on the right-hand side of (A.6) is $o_p(1)$, and the first term $\Lambda_n^{-1/2}\Delta'\psi_\tau(e) \xrightarrow{D} N(0,I_q)$ by assumption A6 and the central limit theorem. $\square$

PROOF OF THEOREM 1. Let

$$\varsigma(\beta,\Theta) = \begin{pmatrix} \varsigma_1 \\ \varsigma_2 \end{pmatrix} = \begin{pmatrix} \Lambda_n^{*-1/2} K_n^*(\beta - \beta_0) \\ k_n^{-1/2} H_n(\Theta - \Theta_0) + k_n^{1/2} H_n^{-1} \Pi' B Z(\beta - \beta_0) \end{pmatrix},$$

$\widehat{\varsigma} = \varsigma(\widehat{\beta}, \widehat{\Theta}) = (\widehat{\varsigma_1}', \widehat{\varsigma_2}')'$, where $H_n^2 = k_n \Pi' B \Pi$, $\beta_0 \in \mathbb{R}^q$ and $\Theta_0 = (\theta_{l0}) \in \mathbb{R}^{p k_n}$. We shall show that $\|\widehat{\varsigma}\| = O_p(k_n^{1/2})$. To do so, we standardize $\widetilde{z}_{ij} = \Lambda_n^{*1/2} K_n^{*-1} z_{ij}^*$, $\widetilde{\pi}_{ij} = k_n^{1/2} H_n^{-1} \Pi_{ij}$. Denote $R_{nij} = \Pi_{ij}' \Theta_0 - \sum_{l=1}^p x_{ij,l} \alpha_l(t_{ij})$ as the bias from the spline approximation. Thus,

$$\sum_{i=1}^n \sum_{j=1}^{m_i} \rho_\tau(y_{ij} - \Pi_{ij}' \Theta - z_{ij}' \beta) = \sum_{ij} \rho_\tau(e_{ij} - \varsigma_1' \widetilde{z}_{ij} - \varsigma_2' \widetilde{\pi}_{ij} - R_{nij}).$$

By Lemma A.5 of [15], $\max_{i,j} \|\widetilde{\pi}_{ij}\| = O(\sqrt{k_n/n})$. Applying similar arguments used in Theorem 3.1 of [25], for any $\epsilon > 0$, there exists $L_\epsilon$ such that

$$P\left\{\inf_{\|\varsigma\| > L_\epsilon k_n^{1/2}} \sum_{ij} \rho_\tau(e_{ij} - \varsigma_1' \widetilde{z}_{ij} - \varsigma_2' \widetilde{\pi}_{ij} - R_{nij}) > \sum_{ij} \rho_\tau(e_{ij} - R_{nij})\right\} > 1 - \epsilon.$$

Using the fact that $\sum_{ij} \rho_\tau(e_{ij} - \varsigma_1' \widetilde{z}_{ij} - \varsigma_2' \widetilde{\pi}_{ij} - R_{nij})$ is minimized at $\hat{\varsigma}$ over the space $R^{p_n}$, we have $P(\|\hat{\varsigma}\| < L_\epsilon k_n^{1/2}) > 1 - \epsilon$, and thus $\|\hat{\varsigma}\| = O_p(k_n^{1/2})$. This together with Lemma 1, and the definition of $\hat{\varsigma}$ gives

$$\|\hat{\beta} - \beta_0\| = \|\Lambda_n^{*1/2} K_n^{*-1} \hat{\varsigma}_1\| = O_p(n^{-1/2} \|\hat{\varsigma}_1\|) = O_p(n^{-1/2} k_n^{1/2}).$$

On the other hand, by (A.1), there exists constants $C_l$, $l = 1, \ldots, p$, such that

$$N^{-1} \sum_{ij} \{\hat{\alpha}_l(t_{ij}) - \alpha_l(t_{ij})\}^2$$

$$\leq 2 N^{-1} \sum_{ij} \{\pi'(t_{ij})(\hat{\theta}_l - \theta_{l0})\}^2 + 2 C_l^2 k_n^{-2r}$$

$$\leq 2 N^{-1} \|\hat{\varsigma}_2\|^2 + 2 \|k_n^{1/2} H_n^{-1} \Pi' B Z (\beta - \beta_0)\|^2 + 2 C_l^2 k_n^{-2r}$$

$$= O_p(n^{-1} \|\hat{\varsigma}_2\|^2) + O_p(\|\hat{\beta} - \beta_0\|^2) + O(k_n^{-2r}) = O_p(k_n^{-2r}).$$

The proof of (2.3) is hence complete.

Next, we will establish the asymptotic normality of $\hat{\beta}$. Let $\hat{\varsigma}_1^* = \Lambda_n^{*-1/2} \times \sum_{i=1}^n Z_i^{*\prime} \psi_\tau(e_i)$. Due to Lemmas 1 and 2, $\hat{\varsigma}_1^*$ is asymptotically normally distributed with variance-covariance matrix $I_q$. Therefore, to show (2.4), all we need is $\|\hat{\varsigma}_1^* - \hat{\varsigma}_1\| = o_p(1)$.

By the definitions of $\hat{\varsigma}_1^*$ and $\hat{\varsigma}_1$, for any $L > 0$ and $M > 0$, we have $P(\hat{\varsigma}_1^* < M) \to 1$ and $P(\hat{\varsigma}_1 < L k_n^{1/2}) \to 1$. Let

$$U_{ij}(\varsigma_1, \hat{\varsigma}_1^*) = \rho_\tau(e_{ij} - \varsigma_1' \widetilde{z}_{ij} - \hat{\varsigma}_2' \widetilde{\pi}_{ij} - R_{nij}) - \rho_\tau(e_{ij} - \hat{\varsigma}_1^{*\prime} \widetilde{z}_{ij} - \hat{\varsigma}_2' \widetilde{\pi}_{ij} - R_{nij}).$$

By Lemmas 8.1 and 8.3 of [25], and the orthogonality of $Z^*$ and $B\Pi$, for any given $\delta > 0$, we have

$$\sup_{\|\varsigma_1 - \hat{\varsigma}_1^*\| < \delta} \left|\sum_{ij} \{U_{ij}(\varsigma_1, \hat{\varsigma}_1^*) + (\varsigma_1 - \hat{\varsigma}_1^*)' \widetilde{z}_{ij} \psi_\tau(e_{ij}) - E U_{ij}(\varsigma_1, \hat{\varsigma}_1^*)\}\right| = o_p(1),$$

$$\sup_{\|\varsigma_1-\hat{\varsigma}_1^*\|<\delta}\left|\sum_{ij}EU_{ij}(\varsigma_1,\hat{\varsigma}_1^*)\right.$$
$$\left.-1/2(\varsigma_1'\Lambda_n^{*1/2}K_n^{*-1}\Lambda_n^{*1/2}\varsigma_1 - \hat{\varsigma}_1^{*\prime}\Lambda_n^{*1/2}K_n^{*-1}\Lambda_n^{*1/2}\hat{\varsigma}_1^*)\right| = o_p(1).$$

Therefore,

$$\sup_{\|\varsigma_1-\hat{\varsigma}_1^*\|<\delta}\left|\sum_{ij}U_{ij}(\varsigma_1,\hat{\varsigma}_1^*) + (\varsigma_1-\hat{\varsigma}_1^*)'\widetilde{Z}'\psi_\tau(e)\right.$$
$$\left.-1/2\varsigma_1'\Lambda_n^{*1/2}K_n^{*-1}\Lambda_n^{*1/2}\varsigma_1 + 1/2\hat{\varsigma}_1^{*\prime}\Lambda_n^{*1/2}K_n^{*-1}\Lambda_n^{*1/2}\hat{\varsigma}_1^*\right|$$

(A.7)
$$= \sup_{\|\varsigma_1-\hat{\varsigma}_1^*\|<\delta}\left|\sum_{ij}U_{ij}(\varsigma_1,\hat{\varsigma}_1^*)\right.$$
$$\left.-1/2(\varsigma_1-\hat{\varsigma}_1^*)'\Lambda_n^{*1/2}K_n^{*-1}\Lambda_n^{*1/2}(\varsigma_1-\hat{\varsigma}_1^*)\right| = o_p(1),$$

where $\widetilde{Z} = (\widetilde{z}_{ij})$. When $\|\varsigma_1 - \hat{\varsigma}_1^*\| > \delta$, $(\varsigma_1 - \hat{\varsigma}_1^*)'(\varsigma_1 - \hat{\varsigma}_1^*) > 0$. Then (A.7) implies that

(A.8)
$$\lim_{n\to\infty}P\left\{\inf_{\|\varsigma_1-\hat{\varsigma}_1^*\|\geq\delta}\sum_{ij}\rho_\tau(e_{ij} - \varsigma_1'\widetilde{z}_{ij} - \hat{\varsigma}_2'\widetilde{\pi}_{ij} - R_{nij})\right.$$
$$\left.> \sum_{ij}\rho_\tau(e_{ij} - \hat{\varsigma}_1^{*\prime}\widetilde{z}_{ij} - \hat{\varsigma}_2'\widetilde{\pi}_{ij} - R_{nij})\right\} = 1.$$

By the convexity of the objective function $\rho_\tau(\cdot)$ and the definition of $\hat{\varsigma}_1$, (A.8) implies that for any $\delta > 0$, $P(\|\hat{\varsigma}_1 - \hat{\varsigma}_1^*\| > \delta) \to 0$, that is, $\|\hat{\varsigma}_1 - \hat{\varsigma}_1^*\| = o_p(1)$. This completes the proof of Theorem 1. □

PROOF OF THEOREM 2. The proof of Theorem 2 follows the similar arguments to those for Theorem 2 of [10], and thus is omitted. □

PROOF OF THEOREM 3. Let $S_n^* = N^{-1/2}\sum_{i=1}^n\{\sum_{j=1}^{m_i}d_{ij}\psi_\tau(e_{ij})\}$, which is a sum of $n$ independent random vectors with mean zero. It's easy to see that

$$\text{Cov}(S_n^*) = N^{-1}\sum_{i=1}^n\text{Cov}\left(\sum_{j=1}^{m_i}d_{ij}\psi_\tau(e_{ij})\right) = N^{-1}\sum_{i=1}^n\mathcal{D}_i'\text{Cov}(\psi_\tau(e_i))\mathcal{D}_i \doteq V_n.$$

It follows from assumption A7 and the central limit theorem that

(A.9) $$V_n^{-1/2}S_n^* \xrightarrow{D} N(0, I_q).$$

By similar arguments used in Theorem 4.1 of [25], we have

$$(A.10) \quad N^{-1/2}\left\|\sum_{ij} \dot{d}_{ij}\{\psi_\tau(y_{ij} - \varpi'_{ij}\hat{\phi}) - \psi_\tau(e_{ij})\}\right\| = o_p(1).$$

Similar to Theorem 2, we can show that under $H_0$, $\hat{V}_n - V_n = o_p(1)$, which together with (A.9), (A.10) completes the proof of Theorem 3. □

PROOF OF THEOREM 4. We first define

$$s_n^* = N^{-1/2}\sum_{i=1}^n \sum_{j=1}^{m_i} d_{ij}\psi_\tau(e_{ij}),$$

where the summands $\sum_{j=1}^{m_i} d_{ij}\psi_\tau(e_{ij})$ are independent of each other and have zero mean. By Theorem 4.1 of Portnoy [20],

$$(A.11) \quad t_n^* = \frac{s_n^{*\prime} v_n^{-1} s_n^* - (k_n + \hbar)p_1}{\sqrt{2(k_n + \hbar)p_1}} \xrightarrow{D} N(0,1),$$

where $v_n = N^{-1}\sum_{i=1}^n D'_i A_i(\Delta) D_i$. Similar to Theorem 2, we can obtain that $\hat{v}_n - v_n = o_p(1)$ under $H_0$. Due to (A.11), it is clear that all we need to show is

$$(A.12) \quad N^{-1/2}\left\|\sum_{i=1}^n \sum_{j=1}^{m_i} d_{ij}\{\psi_\tau(\hat{e}_{ij}) - \psi_\tau(e_{ij})\}\right\| = o_p(1).$$

Let $\varphi = (\xi'_2, \beta')'$, and $\varphi_0 = (\xi'_{20}, \beta'_0)'$ such that $y_{ij} = \Pi_{ij}^{(1)\prime}\xi_{10} + w'_{ij}\varphi_0 + e_{ij} + R_{nij}$, where $R_{nij}$ is the bias from the B-spline approximation satisfying $R_{nij} = O(k_n^{-r})$. Denote

$$u_{ij}(\varphi, \varphi_0) = d_{ij}[\psi_\tau(y_{ij} - w'_{ij}\varphi) - \psi_\tau(y_{ij} - w'_{ij}\varphi_0) \\ - E\{\psi_\tau(y_{ij} - w'_{ij}\varphi) - \psi_\tau(y_{ij} - w'_{ij}\varphi_0)\}].$$

In our context, the dimensions of both $d_{ij}$ and $w_{ij}$ grow with $n$, while the former is finite in [25], and the latter is finite in [15]. By Assumption A.4 and Lemma A.5 of [15],

$$(A.13) \quad \sup_{ij}\|d_{ij}\| = O_p(\sqrt{k_n}), \quad \sup_{ij}\|w_{ij}\| = O_p(\sqrt{k_n}).$$

Let $\delta_n$ be a sequence of positive numbers. It follows from (A.13), and Lemmas 2.2 and 3.3 of He and Shao [9] that for any $L > 0$,

$$(A.14) \quad \sup_{\|\varphi - \varphi_0\| \leq L\delta_n} N^{-1/2}\left\|\sum_{ij} u_{ij}(\varphi, \varphi_0)\right\| = O_p(k_n(\delta_n \log n)^{1/2}).$$

By expanding $E\{d_{ij}\psi_\tau(y_{ij} - w'_{ij}\varphi)\}$ around $\varphi_0$ for each $i, j$, and using (A.13) and the fact that $D'BW = 0$, we have

$$
\sup_{\|\varphi-\varphi_0\|\leq L\delta_n} N^{-1/2}\left\|\sum_{ij} E[d_{ij}\{\psi_\tau(y_{ij} - w'_{ij}\varphi) - \psi_\tau(y_{ij} - w'_{ij}\varphi_0)\}]\right\|
$$
(A.15)
$$
= O(n^{1/2}(k_n^{-r+1}\delta_n + k_n^{3/2}\delta_n^2)).
$$

Theorem 1 provides the consistency of $\hat{\varphi}$ for $k_n \approx n^{1/(2r+1)}$, which, however, is not the necessary condition for the consistency. In fact, for any $k_n = o(n^{1/3})$, we have $\|\hat{\varphi} - \varphi_0\| = O_p(n^{-1/2}k_n^{1/2})$. Therefore, it follows from (A.14) and (A.15) that

(A.16) $\quad N^{-1/2}\left\|\sum_{ij} d_{ij}\{\psi_\tau(y_{ij} - w'_{ij}\hat{\varphi}) - \psi_\tau(y_{ij} - w'_{ij}\varphi_0)\}\right\| = o_p(1).$

To show (A.12), it suffices to show

(A.17) $\quad N^{-1/2}\sum_{i=1}^{n}\sum_{j=1}^{m_i}\|d_{ij}\{\psi_\tau(y_{ij} - w'_{ij}\varphi_0) - \psi_\tau(e_{ij})\}\| = o_p(1).$

Note that under $H_0$, $y_{ij} - w'_{ij}\varphi_0 = e_{ij} + R_{nij}$. Therefore, the uniform approximation technique used to obtain (A.14) is not applicable here, as the left-hand side of (A.17) involves unknown bias with dimension of the same order as $n$. An alternative way to prove (A.17) is to show the $L_1$ convergence of the left-hand side of (A.17), which, however, is difficult for $k_n \ll n^{1/(2r-1)}$, including $k_n \approx n^{1/(2r+1)}$; see (A.23). To bypass the technical difficulty, we take an intermediate step.

Let $\aleph$ be the set of knots used in the estimation. Denote $k_n$ as the dimension of $\aleph$. For easy demonstration, we assume that $k_n \approx n^{1/(2r+1)}, r > 2$, but the same proof goes through for other $k_n$. Similar to (A.1), we obtain that under $H_0$, $\sup_{ij}|R_{nij}| = O(k_n^{-r})$. By adding more knots into $\aleph$, we have a new set of knots $\widetilde{\aleph}$ with its dimension $\widetilde{k}_n \approx n^\alpha$, where $n^{(r+1)/\{(2r+1)r\}} \ll \widetilde{k}_n \ll n^{1/5}$. As $\aleph$ is a subsequence of $\widetilde{\aleph}$, the original B-spline space is a subspace of the new one, for which we can construct a set of basis functions such that the basis functions for the original space is a subset and is orthogonal to the additional basis functions. Using the new set of knots, we define $\widetilde{W}, \widetilde{w}_{ij}, \widetilde{\varphi}, \widetilde{\varphi}_0$ and $\widetilde{R}_{nij}$ the same way as $W, w_{ij}, \hat{\varphi}, \varphi_0$ and $R_{nij}$. Let $\widetilde{P}_w = \widetilde{W}(\widetilde{W}'B\widetilde{W})^{-1}\widetilde{W}'B$, $\widetilde{D} = (I - \widetilde{P}_w)\Pi^{(1)}$, and $\widetilde{d}_{ij}$ be the row components of $\widetilde{D}$ corresponding to the $j$th measurement of the $i$th subject. As the $(\hbar+1)$th order B-splines are $(\hbar-1)$-times differentiable, similar to (A.1), for $\hbar \geq 3$, we have

(A.18) $\quad\sup_{ij}\|d_{ij} - \widetilde{d}_{ij}\| = O_p(k_n^{-\hbar+2} + \widetilde{k}_n^{-\hbar+1}k_n).$

Using (A.18), and the facts that $E\{\psi_\tau(y_{ij} - w'_{ij}\varphi_0)\} = E\{\psi_\tau(R_{ij} + e_{ij})\} = O(R_{nij}) = O(k_n^{-r})$ and $E\{\psi_\tau(y_{ij} - \widetilde{w}'_{ij}\widetilde{\varphi}_0)\} = O(\widetilde{k}_n^{-r})$, we have

$$N^{-1/2}\left\|\sum_{ij}(d_{ij} - \widetilde{d}_{ij})\psi_\tau(y_{ij} - w'_{ij}\varphi_0)\right\| = o_p(1) \quad \text{and}$$

(A.19)

$$N^{-1/2}\left\|\sum_{ij}(d_{ij} - \widetilde{d}_{ij})\psi_\tau(y_{ij} - \widetilde{w}'_{ij}\widetilde{\varphi}_0)\right\| = o_p(1).$$

Note that $\sup_{ij}|w'_{ij}\varphi_0 - \widetilde{w}'_{ij}\widetilde{\varphi}_0| \leq \sup_{ij}|R_{nij}| + \sup_{ij}|\widetilde{R}_{nij}| = O(k_n^{-r})$. Applying the similar arguments used for (A.14), we can show that

$$N^{-1/2}\left\|\sum_{ij}\widetilde{d}_{ij}\{\psi_\tau(y_{ij} - w'_{ij}\varphi_0) - \psi_\tau(y_{ij} - \widetilde{w}'_{ij}\widetilde{\varphi}_0)\}\right.$$

(A.20)

$$\left. - \sum_{ij}E[\widetilde{d}_{ij}\{\psi_\tau(y_{ij} - w'_{ij}\varphi_0) - \psi_\tau(y_{ij} - \widetilde{w}'_{ij}\widetilde{\varphi}_0)\}]\right\| = o_p(1).$$

For each $i$ and $j$, we expand the conditional mean $E\{\widetilde{d}_{ij}\psi_\tau(y_{ij} - \widetilde{w}'_{ij}\widetilde{\varphi}_0)|(\widetilde{d}_{ij}, \widetilde{w}_{ij}, w_{ij})\}$ around $w'_{ij}\varphi_0$. Recall that $\sup_{i,j}\|\widetilde{d}_{ij}\| = O_p(\widetilde{k}_n^{1/2})$, $\sup_{ij}|w'_{ij}\varphi_0 - \widetilde{w}'_{ij}\widetilde{\varphi}_0| = O_p(k_n^{-r})$, $y_{ij} - w'_{ij}\varphi_0 = e_{ij} + R_{nij}$ with $\sup_{ij}|R_{nij}| = O_p(k_n^{-r})$. Under assumption A4, we have

$$E[\widetilde{d}_{ij}\{\psi_\tau(y_{ij} - w'_{ij}\varphi_0) - \psi_\tau(y_{ij} - \widetilde{w}'_{ij}\widetilde{\varphi}_0)\}|(\widetilde{d}_{ij}, w_{ij}, \widetilde{w}_{ij})]$$
$$= \widetilde{d}_{ij}f_{ij}(R_{nij})(\widetilde{w}'_{ij}\widetilde{\varphi}_0 - w'_{ij}\varphi_0) + O_p(\widetilde{k}_n^{1/2}k_n^{-2r})$$
$$= \widetilde{d}_{ij}f_{ij}(0)(\widetilde{w}'_{ij}\widetilde{\varphi}_0 - w'_{ij}\varphi_0) + O_p(\widetilde{k}_n^{1/2}k_n^{-2r}) \quad \text{uniformly over } i \text{ and } j.$$

Therefore,

$$N^{-1/2}\left\|\sum_{ij}E[\widetilde{d}_{ij}\{\psi_\tau(y_{ij} - w'_{ij}\varphi_0) - \psi_\tau(y_{ij} - \widetilde{w}'_{ij}\widetilde{\varphi}_0)\}]\right\|$$

(A.21)
$$= N^{-1/2}\|E(\widetilde{D}'B\widetilde{W}\widetilde{\varphi}_0) - E(\widetilde{D}'BW\varphi_0)\| + O(n^{1/2}\widetilde{k}_n^{1/2}k_n^{-2r})$$
$$= o(1),$$

where the last step comes from the facts that $\widetilde{D}$ is orthogonal to both $\widetilde{W}'B$ and $W'B$, and that $n^{1/2}\widetilde{k}_n^{1/2}k_n^{-2r} = o(1)$ for $k_n = n^{1/(2r+1)}$ with $r > 2$ and $\widetilde{k}_n \ll n^{1/5}$. Collecting (A.16) and (A.19)–(A.21), we get

(A.22) $\quad N^{-1/2}\left\|\sum_{ij}d_{ij}\{\psi_\tau(y_{ij} - w'_{ij}\hat{\varphi}) - \psi_\tau(y_{ij} - \widetilde{w}'_{ij}\widetilde{\varphi}_0)\}\right\| = o_p(1).$

...

Furthermore, note that $|\psi_\tau(e_{ij} + \widetilde{R}_{nij}) - \psi_\tau(e_{ij})| \leq I(|e_{ij}| \leq |\widetilde{R}_{nij}|)$. By (A.18), we have for $\widetilde{k}_n \gg n^{(r+1)/\{(2r+1)r\}}$,

$$N^{-1/2} E \sum_{ij} \|d_{ij}\{\psi_\tau(y_{ij} - \widetilde{w}'_{ij}\widetilde{\varphi}_0) - \psi_\tau(e_{ij})\}\|$$

(A.23)
$$\leq N^{-1/2} E \sum_{ij} |\widetilde{R}_{nij}| \cdot \|d_{ij}\| = O(\widetilde{k}_n^{-r}(k_n n)^{1/2}) = o(1),$$

which together with (A.22) gives (A.12). The proof of Theorem 4 is thus complete. $\square$

**Acknowledgments.** The authors are grateful to two referees and an Associate Editor for their helpful comments and suggestions that led to an improvement of this article.
## REFERENCES

[1] Ahmad, I., Leelahanon, S. and Li, Q. (2005). Efficient estimation of a semiparametric partially linear varying coefficient model. *Ann. Statist.* **33** 258–283. MR2157803
[2] Cai, Z. and Xu, X. (2008). Nonparametric quantile estimations for dynamic smooth coefficient models. *J. Amer. Statist. Assoc.* **103** 1595–1608.
[3] Chiang, C. T., Rice, J. A. and Wu, C. O. (2001). Smoothing spline estimation for varying-coefficient models with repeatedly measured dependent variables. *J. Amer. Statist. Assoc.* **96** 605–619. MR1946428
[4] Fan, J. and Huang, T. (2005). Profile likelihood inferences on semiparametric varying-coefficient partially linear models. *Bernoulli* **11** 1031–1057. MR2189080
[5] Fan, J., Huang, T. and Li, R. (2007). Analysis of longitudinal data with semiparametric estimation of covariance function. *J. Amer. Statist. Assoc.* **102** 632–641. MR2370857
[6] Gutenbrunner, C., Jurêcková, J., Koenker, R. and Portnoy, S. (1993). Tests of linear hypotheses based on regression rank scores. *J. Nonparametr. Stat.* **2** 307–333. MR1256383
[7] Hall, P. and Sheather, S. J. (1988). On the distribution of a studentized quantile. *J. Roy. Statist. Soc. Ser. B* **50** 381–391. MR0970974
[8] Hastie, T. and Tibshirani, R. (1993). Varying-coefficient models. *J. Roy. Statist. Soc. Ser. B* **55** 757–796. MR1229881
[9] He, X. and Shao, Q. M. (2000). On parameters of increasing dimensions. *J. Multivariate Anal.* **73** 120–135. MR1766124
[10] He, X., Zhu, Z. Y. and Fung, W. K. (2002). Estimation in a semiparametric model for longitudinal data with unspecified dependence structure. *Biometrika* **89** 579–590. MR1929164
[11] Hendricks, W. and Koenker, R. (1992). Hierarchical spline models for conditional quantiles and the demand for electricity. *J. Amer. Statist. Assoc.* **87** 58–68.
[12] Honda, T. (2004). Quantile regression in varying coefficient models. *J. Statist. Plann. Inference* **121** 113–125. MR2027718
[13] Huang, J. Z., Wu, C. O. and Zhou, L. (2002). Varying-coefficient models and basis function approximations for the analysis of repeated measurements. *Biometrika* **89** 111–128. MR1888349

[14] KASLOW, R. A., OSTROW, D. G., DETELS, R., PHAIR, J. P., POLK, B. F. and RINALDO, C. R. (1987). The multicenter AIDS cohort study: Rationale, organization and selected characteristics of the participants. *American Journal of Epidemiology* **126** 310–318.

[15] KIM, M. (2007). Quantile regression with varying coefficients. *Ann. Statist.* **35** 92–108. MR2332270

[16] KOENKER, R. (1994). Confidence intervals for regression quantiles. In *Asymptotic Statistics: Proceedings of the 5th Prague Symposium on Asymptotic Statistics* (P. Mandl and M. Husková, eds.) 349–359. Physica, Heidelberg. MR1311953

[17] KOENKER, R. (2004). Quantile regression for longitudinal data. *J. Multivariate Anal.* **91** 74–89. MR2083905

[18] LIPSITZ, S. R., FITZMAURICE, G. M., MOLENBERGHS, G. and ZHAO, L. P. (1997). Quantile regression methods for longitudinal data with drop-outs: Application to CD4 cell counts of patients infected with the human immunodeficiency virus. *J. Roy. Statist. Soc. Ser. C* **46** 463–476.

[19] MU, Y. and WEI, Y. (2009). A dynamic quantile regression transformation model for longitudinal data. *Statist. Sinica* **19** 1137–1153.

[20] PORTOY, S. (1985). Asymptotic behavior of $M$-estimators of $p$ regression parameters when $p^2/n$ is large; II. Normal approximation. *Ann. Statist.* **13** 1403–1417. MR0811499

[21] QU, A. and LI, R. (2006). Quadratic inference functions for varying coefficient models with longitudinal data. *Biometrics* **62** 379–391. MR2227487

[22] SCHUMAKER, L. L. (1981). *Spline Functions: Basic Theory*. Wiley, New York. MR0606200

[23] SUN, Y. and WU, H. (2005). Semiparametric time-varying coefficients regression model for longitudinal data. *Scand. J. Statist.* **32** 21–47. MR2136800

[24] WANG, H. and HE, X. (2007). Detecting differential expressions in GeneChip microarray studies: A quantile approach. *J. Amer. Statist. Assoc.* **102** 104–112. MR2293303

[25] WEI, Y. and HE, X. (2006). Conditional growth charts (with discussion). *Ann. Statist.* **34** 2069–2097. MR2291494

[26] YU, K. and JONES, M. C. (1998). Local linear quantile regression. *J. Amer. Statist. Assoc.* **93** 228–237. MR1614628

H. J. WANG
DEPARTMENT OF STATISTICS
NORTH CAROLINA STATE UNIVERSITY
RALEIGH, NC 27695
USA
E-MAIL: wang@stat.ncsu.edu

Z. ZHU
DEPARTMENT OF STATISTICS
FUDAN UNIVERSITY
SHANGHAI, 200433
CHINA
E-MAIL: zhuzy@fudan.edu.cn

J. ZHOU
DEPARTMENT OF STATISTICS
UNIVERSITY OF VIRGINIA
CHARLOTTESVILLE, VA 22904
E-MAIL: jz9p@virginia.edu